\documentclass[graybox]{svmult}

\usepackage{type1cm}        
\usepackage{makeidx}         
\usepackage{graphicx}        
\usepackage{multicol}        
\usepackage[bottom]{footmisc}

\usepackage{newtxtext}       %
\usepackage[varvw]{newtxmath}       

\makeindex             

\newcommand{\bbE}{{\mathbb{E}}}
\newcommand{\bbV} {{\mathbb{V}}}
\newcommand{\bbP}{{\mathbb{P}}}
\newcommand{\bsone}{{\boldsymbol{1}}}

\graphicspath{{figs/}}

\begin{document}

\title*{MLMC techniques for discontinuous functions}
\author{Michael B.~Giles}
\institute{Michael B.~Giles \at University of Oxford Mathematical Institute,
           Woodstock Rd, Oxford OX2 6GG, UK\\
  \email{mike.giles@maths.ox.ac.uk}}
%
%
\maketitle

\abstract*{The Multilevel Monte Carlo (MLMC) approach usually works well
           when estimating the expected value of a quantity which is a
           Lipschitz function of intermediate quantities, but if it is
           a discontinuous function it can lead to a much slower decay in
           the variance of the MLMC correction.
           This article reviews the literature on techniques which can
           be used to overcome this challenge in a variety of different
           contexts, and discusses recent developments using either a 
           branching diffusion or adaptive sampling.}

\abstract{The Multilevel Monte Carlo (MLMC) approach usually works well
           when estimating the expected value of a quantity which is a
           Lipschitz function of intermediate quantities, but if it is
           a discontinuous function it can lead to a much slower decay in
           the variance of the MLMC correction.
           This article reviews the literature on techniques which can
           be used to overcome this challenge in a variety of different
           contexts, and discusses recent developments using either a 
           branching diffusion or adaptive sampling.}

\section{Introduction}

The Multilevel Monte Carlo (MLMC) method is based on the telescoping sum

\[
\bbE[{\widehat{P}}_L] = \bbE[{\widehat{P}}_0] + \sum_{\ell=1}^L \bbE[{\widehat{P}}_\ell{-}{\widehat{P}}_{\ell-1}]
\]
where ${\widehat{P}}_\ell$ represents an approximation to an output quantity
of interest $P$ on level $\ell$, with the weak error
$\left|\bbE[{\widehat{P}}_\ell{-}P]\right|$ and MLMC variance $\bbV[{\widehat{P}}_\ell{-}{\widehat{P}}_{\ell-1}]$,
both decreasing as the level $\ell$ increases, but with the corresponding
computational cost per sample increasing.

If ${\widehat{Y}}_\ell$ has expected value $\bbE[{\widehat{P}}_\ell{-}{\widehat{P}}_{\ell-1}]$,
with ${\widehat{P}}_{-1}{\equiv} 0$, 
with variance $V_\ell$ and cost $C_\ell$, then for a given target RMS error
$\varepsilon$, the number of levels $L$ and the number of independent samples
on each level can be optimised \cite{Gi_giles08,Gi_giles15} to give an 
overall cost which is approximately equal to 
$\displaystyle
  \varepsilon^{-2} \left( \sum_{\ell=0}^L \sqrt{C_\ell V_\ell} \right)^2.
$

If $C_\ell V_\ell \rightarrow 0$ as $\ell\rightarrow \infty$, then the cost
is dominated by the first term from level $0$, and the cost is approximately
$\varepsilon^{-2} C_0 V_0$, so proportional to $\varepsilon^{-2}$.

If $C_\ell V_\ell \rightarrow \mbox{\rm const}$ as
$\ell\rightarrow \infty$, then the contributions to the cost are
spread almost equally across all levels and the cost is approximately
$\varepsilon^{-2} L^2 C_L V_L$, proportional to $\varepsilon^{-2}|\log\varepsilon|^2$
if $\bbE[{\widehat{P}}_\ell{-}P]$ decays exponentially with $\ell$.

Even worse, if $C_\ell V_\ell \rightarrow \infty$ as $\ell\rightarrow \infty$,
then the cost is dominated by the contribution from the finest level and so is
approximately $\varepsilon^{-2} C_L V_L$ which is $O(\varepsilon^{-2 - (\gamma-\beta)/\alpha})$
if $\bbE[{\widehat{P}}_\ell{-}P] \sim 2^{-\alpha \ell}$, $V_\ell \sim 2^{-\beta \ell}$ and $C_\ell \sim 2^{\gamma \ell}$.

In most MLMC applications, $P$ is a smooth, or at least Lipschitz,
function of some intermediate solution quantities, such as the solution
of an SDE, a PDE with stochastic coefficients or initial/boundary data,
or an estimate of an inner conditional expectation.  Under these
circumstances we usually have $\beta \geq \gamma$ and so the
MLMC complexity is $O(\varepsilon^{-2})$ or $O(\varepsilon^{-2}|\log\varepsilon|^2)$.

This article is concerned with the small but important class of
applications where $P$ is a discontinuous function of the intermediate
quantities, and because of this the MLMC variance $V_\ell$ can decay
much more slowly, leading to the complexity falling into the third
category of being $O(\varepsilon^{-2 - (\gamma-\beta)/\alpha})$.

The good news is that there has been considerable research within
the MLMC community
\footnote{See \tt https://people.maths.ox.ac.uk/gilesm/mlmc\_community.html}
to address this challenge.  This article surveys
the variety of methods which have been developed in the hope that this
can aid researchers meeting similar challenges in future applications.

To illustrate things, we begin by detailing two specific application
challenges which have motivated much of this research.  We then discuss
the many approaches which have been developed, several of which have
borrowed ideas from the literature on computing sensitivities
(the ``greeks'' in mathematical finance literature) of the form
$
\frac{\partial }{\partial \alpha} \bbE\left[f(\omega,\alpha)\right]
$
using the pathwise sensitivity approach \cite{Gi_glasserman04}
(also known as Infinitesimal Perturbation Analysis, IPA for short
\cite{Gi_ecuyer90}) or Likelihood Ratio Method \cite{Gi_ecuyer95},
or from methods for improving integrand smoothness to improve
the rate of convergence for QMC integration \cite{Gi_acn13,Gi_bst18}.

\subsection{Challenge 1: nested expectation}

Suppose $f$ is a scalar function and we want to estimate the nested expectation 
$\displaystyle \bbE \left[ f\left(\bbE[Z | X]\right)  \, \right] $,
where the outer expectation is with respect to a random variable $X$
and we will assume that the inner conditional expectation $\bbE[Z | X]$
has a bounded density near zero.

A very simple MLMC treatment
\footnote{Note that if $f$ is smooth, or at least Lipschitz, then it is
better to use an ``antithetic'' estimator \cite{Gi_bhr15,Gi_giles15,Gi_giles18,Gi_gg19},
but this does not give a better order of convergence when $f$ is discontinuous.}
uses $M_\ell=2^\ell M_0$ inner samples on level $\ell$, so estimators
on level $0$ and the higher levels are simply
\[
{\widehat{Y}}_0 = f( {\overline{Z}}^{(0,M_0)}), ~~~~~
{\widehat{Y}}_\ell = f( {\overline{Z}}^{(\ell,M_\ell)}) - f( {\overline{Z}}^{(\ell,M_{\ell-1})}),
\]
where $ {\overline{Z}}^{(\ell,M_\ell)}$ and $ {\overline{Z}}^{(\ell,M_{\ell-1})}$ represent
independent averages of $M_\ell$ and $M_{\ell-1}$ independent
samples of $Z$, all conditional on the same value of $X$
\cite{Gi_giles15,Gi_giles18}.

If $\bbV[Z | X]$ is finite, and $f$ is Lipschitz with constant $L_f$, then
\begin{eqnarray*}
\bbE\left[\left( f( {\overline{Z}}^{(\ell,M_\ell)})
                    - f\left( \bbE[ Z | X] \right)  \right)^2 |\ X \right]
 &\leq& L_f^2\  \bbE\left[\left(
      {\overline{Z}}^{(\ell,M_\ell)} - \bbE[ Z | X] \right)^2 | \ X \right] \\
 &=& L_f^2 M_\ell^{-1} \bbV[  Z | X ],
\end{eqnarray*}
and hence
$\bbE[{\widehat{Y}}^2_\ell | X]\leq 2\,L_f^2 (M_\ell^{-1}+M^{-1}_{\ell-1}) \bbV[  Z | X ]$
for $\ell{>}0$.
If $\bbV[  Z | X ]$ is uniformly bounded it follows that
$V_\ell = O(M_\ell^{-1})$. If the cost of each conditional
sample of $Z$ is $O(1)$ then $C_\ell=O(M_\ell)$ and hence
the complexity is $O(\varepsilon^{-2}|\log \varepsilon|^2)$.

Unfortunately, the situation is significantly poorer when
$f$ is the Heaviside step function $H$ defined by
$H(x){=}0$ if $x{<}0$, and $H(x){=}1$ if $x{\geq} 0$.
This occurs in many applications because
$
\bbP\left[
  \bbE[Z|X] > K \right] =
\bbE\left[
  H(\bbE[Z|X]-K)\right],
$
so it corresponds to the probability of a conditional expectation
exceeding some threshold $K$, which is a very important quantity
in risk calculations.

If $K{=}0$ and $E[Z|X]$ has a bounded density
near zero then there is an $O(M_\ell^{-1/2})$ probability that
$|E[Z|X]\,| = O(M_\ell^{-1/2})$, which is the circumstance under
which there is an $O(1)$ probability that ${\widehat{Y}}_\ell = \pm 1$
due to $ {\overline{Z}}^{(\ell,M_\ell)}$ being positive and $ {\overline{Z}}^{(\ell,M_{\ell-1})}$
being negative, or vice versa. Hence $V_\ell \approx O(M_\ell^{-1/2})$
and the complexity is approximately $O(\varepsilon^{-5/2})$ \cite{Gi_gh19}.

This challenge is the primary motivation for Section 7, and also
arises in the context of Section 3.

\subsection{Challenge 2: discontinuous payoff function}

In the case of a scalar SDE
\begin{equation}
  \D S_t = a(S_t) \, \D t + b(S_t)\, \D W_t,
  \label{Gi_eq:SDE}
\end{equation}
with an output quantity of interest $P\equiv f(S_T)$,
the standard estimator is
\[
{\widehat{Y}}_\ell = {\widehat{P}}_\ell - {\widehat{P}}_{\ell-1}
\]
where the same Brownian motion $W_t$ is used to calculate both
${\widehat{P}}_\ell$ and ${\widehat{P}}_{\ell-1}$, but with different uniform timesteps
$h_\ell$ and $h_{\ell-1}$.

If $f$ is Lipschitz with constant $L_f$, then
\[
  V_\ell \ \leq\ \bbE\left[ ({\widehat{P}}_\ell - {\widehat{P}}_{\ell-1})^2\right]
  \ \leq\ L_f^2\ \bbE\left[ ({\widehat{S}}_\ell - {\widehat{S}}_{\ell-1})^2\right]
\]
where ${\widehat{S}}_\ell$ is the level $\ell$ numerical approximation to $S_T$.
Hence, based on standard strong convergence results \cite{Gi_kp92}
we have $V_\ell=O(h_\ell)$ for an Euler-Maruyama discretisation of the SDE, 
and $V_\ell=O(h_\ell^2)$ for the first order Milstein discretisation.
The cost $C_\ell$ is $O(h_\ell^{-1})$ in both cases, giving MLMC complexities
of $O(\varepsilon^{-2}|\log\varepsilon|^2)$ and $O(\varepsilon^{-2})$, respectively,

In mathematical finance, a digital call option payoff is $0$ or $1$,
depending on whether $S_T$ is below or above the strike $K$,
so the payoff function can be written as $f(S_T)=H(S_T{-}K)$.
The MLMC problem is that a small difference between the coarse and
fine paths can give a payoff difference of $\pm 1$ if the two paths
straddle the strike, i.e.~are on different sides of the strike.

When using the Euler-Maruyama approximation of the SDE,
${\widehat{S}}_\ell{-}{\widehat{S}}_{\ell-1} = O(h_\ell^{1/2})$.
Speaking loosely (see \cite{Gi_avikainen09,Gi_ghm09} for the rigorous analysis)
an $O(h_\ell^{1/2})$ fraction of fine/coarse pairs straddle the strike,
so $V_\ell = O(h_\ell^{1/2})$ and hence the complexity is $O(\varepsilon^{-5/2})$.

Similarly, using the Milstein approximation gives
${\widehat{S}}_\ell{-}{\widehat{S}}_{\ell-1} = O(h_\ell)$
so $V_\ell = O(h_\ell)$. This is clearly better, and gives a complexity
which is $O(\varepsilon^{-2}|\log\varepsilon|^2)$, but there is still the
problem that most MLMC samples $Y_\ell$ are zero on the finer levels,
so the kurtosis is $O(h_\ell^{-1})$ which causes problems in practice
in estimating $V_\ell$ accurately to determine the number of samples
$N_\ell$ to use on level $\ell$.  In addition, there is the difficulty
that the Milstein discretisation of multi-dimensional SDEs often
requires the simulation of L\'evy areas, though this problem can be
addressed through the use of an antithetic estimator \cite{Gi_gs14}.

This challenge is the primary motivation for Sections 2, 4, 5 and 6,
also also arises in Sections 3 and 7.

\section{Explicit smoothing}

The pathwise sensitivity analysis (or IPA) approach to compute the
parameter sensitivities known as Greeks in mathematical finance
\cite{Gi_glasserman04} requires that the payoff function $f$ is continuous
and piecewise smooth. This is clearly a problem with digital options,
and one standard approach in this case is to smooth the payoff function by replacing
the Heaviside step function $H$ with a smoothed approximation
$H_\delta(x) \equiv g(x/\delta)$, with $g(x) \rightarrow 0$ as
$x\rightarrow -\infty$ and $g(x) \rightarrow 1$ as $x\rightarrow +\infty$,
so the discontinuity is smoothed over a width of $O(\delta)$.

\begin{figure}[b]
\begin{center}
\includegraphics[width=0.8\textwidth]{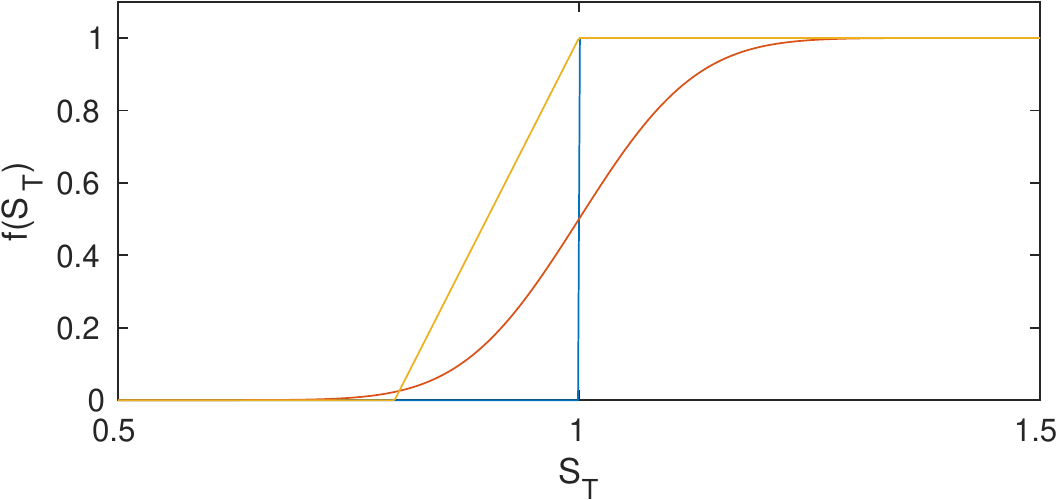}
\end{center}
\caption{Two explicitly smoothed versions of the Heaviside step
  function for a digital call option}
\label{Gi_fig:smoothing}
\end{figure}

For financial reasons, the preference is often to use a one-sided smoothing,
such as the piecewise linear approximation shown in yellow in Figure
\ref{Gi_fig:smoothing}.
This one-sided approximation introduces a weak error, or bias, which
is $O(\delta)$.  If it is used for MLMC, then $H'_\delta(S_T)=\delta^{-1}$
for the $O(\delta)$ fraction of the paths which end up in the ramp
region, and therefore
$V_\ell = O(\delta \times (\delta^{-1})^2) = O(\delta^{-1})$.
Hence the optimal choice of $\delta$ involves a tradeoff
between bias and variance.

The bias can be reduced by making the smoothing anti-symmetric
about $x=0$ so that $H_\delta(x)-H(x) = - (H_\delta(-x)-H(-x))$,
for example by choosing $g(x)\equiv \Phi(x)$ as illustrated in orange
in Figure \ref{Gi_fig:smoothing}.
If $S_T$ has the smooth probability density $\rho(S)$ then the
weak error is
\[
\int_{-\infty}^\infty \left(H_\delta(S{-}K)-H(S{-}K)\right) \rho(S) \, \D S
  =
  \delta \int_{-\infty}^\infty \left( g(x)-H(x) \right)\, \rho(K {+} x \delta) \, \D x
  \]
and a Taylor series expansion of $\rho(K {+} x \delta)$ results in the
asymptotic error expansion
\[
a_1 \rho(K) \, \delta +
a_2 \rho'(K) \, \delta^2 +
a_3 \rho''(K) \, \delta^3 +
a_4 \rho'''(K) \, \delta^4 + O(\delta^5)
\]
where
\[
a_k
= \int_{-\infty}^\infty x^{k-1}  \left( g(x)-H(x) \right) \, \D x.
\]
If $g(x)-H(x) = -\, (g(-x)-H(-x))$ then $a_1=a_3=0$, and
\[
a_2 = 2 \int_0^\infty x (g(x)-1) \, \D x, ~~~
a_4 = 2 \int_0^\infty x^3 (g(x)-1) \, \D x.
\]
If $g(x)$ is monotonic, then $a_2\neq 0$, but by considering
non-monotonic functions such as
$g(x) = (4/3)\, \Phi(x) - (1/3)\, \Phi(2x)$
it is possible to set $a_2=0$ making the weak error $O(\delta^4)$.
Hence, to achieve $O(\varepsilon)$ accuracy overall we need
$\delta {=} O(\varepsilon^{1/4})$, and then on the coarsest levels
$V_\ell = O(\delta^{-1}) = O(\varepsilon^{-1/4})$
so the overall complexity is approximately $O(\varepsilon^{-2-1/4})$ in
the best cases where the overall cost is dominated by the cost
on the coarsest levels.

Giles, Nagapetyan \& Ritter \cite{Gi_gnr15} used explicit smoothing
for estimating cumulative distribution functions (CDFs). For a scalar
random variable $X$, to estimate $C(x) = \bbP(X<x) = \bbE[ H(x{-}X)]$,
the approach they adopted was to use MLMC to estimate $C(x_j)$ for
a set of spline points $x_j$ and then interpolate these values with
a cubic spline.  Overall, their method balanced three weak errors,
the SDE discretisation error on the finest level, the smoothing
error due to $H_\delta$, and the cubic spline interpolation error,
in addition to the MLMC sampling error.


\section{Integration/differentiation and Malliavin calculus}

Krumscheid \& Nobile \cite{Gi_kn18} used a slightly different approach
for estimating CDFs, particularly in the context of risk estimation.
Starting from the identity
\[
\frac{\D}{\D x}\ \bbE\left[\, \max(0,x{-}S_T)\, \right] =  \bbE[\, H(x{-}S_T) \,] 
\]
they used MLMC to estimate $\bbE[\, \max(0,x_j{-}S_T) \,]$ for a set of
spline points $x_j$, interpolated these with a cubic spline, and 
and then differentiated the spline to obtain an approximation to the
CDF $C(x)$.
This avoids the extra weak error due to smoothing the Heaviside function,
but differentiating the cubic spline amplifies the noise in the spline data.

On a similar note, Altmayer \& Neuenkirch \cite{Gi_an15} used Malliavin
calculus integration by parts to treat discontinuous payoffs based on
solutions of the Heston stochastic volatility SDE.  They observed that
asymptotically this improves the MLMC variance on the finer levels, but
it increases the variance on coarse levels.
To address this, they split the payoff into a smooth part which they treated
with the standard MLMC approach, and a compact-support discontinuous part
for which they used the Malliavin MLMC.

Malliavin calculus was originally developed for computing sensitivities,
so this is another example of the literature on sensitivity calculations
being exploited to develop improved MLMC algorithms.

\section{Conditional expectation}

When using the first order Milstein discretisation for an SDE, one way
to improve the MLMC variance for digital options is to switch to the
Euler-Maruyama approximation for the final timestep, and then take the
conditional expectation with respect to the final fine path Brownian
increment $\Delta W$ \cite{Gi_giles08b,Gi_gdr19}.

For the fine path approximation of the scalar SDE (\ref{Gi_eq:SDE}) with
$N$ timesteps of size $h_\ell$, the path value $S_T$ at the final time
$T$ is given by
\[
{\widehat{S}}^f_T = {\widehat{S}}^f_{T-h_\ell} + a({\widehat{S}}^f_{T-h_\ell})\, h_\ell +  b({\widehat{S}}^f_{T-h_\ell})\, \Delta W_N,
\]
and therefore the conditional expected value for the digital call option
is
\[
{\widehat{P}}^f_\ell
\ =\ \bbE\left[ H({\widehat{S}}^f_T {-} K)\ |\ {\widehat{S}}^f_{T-h_\ell}\right]
\ =\ \Phi\left(\frac{{\widehat{S}}^f_{T-h_\ell} + a({\widehat{S}}^f_{T-h_\ell})\, h_\ell - K}
                    {b({\widehat{S}}^f_{T-h_\ell})\ \sqrt{h_\ell}}\right).
\]
Similarly, for the coarse path with coarse timestep $h_{\ell-1}=2\,h_\ell$,
the Brownian increment for the final coarse timestep is the sum of the
last two Brownian increments for the fine path, $\Delta W_{N-1}{+}\Delta W_N$,
and therefore
\[
  {\widehat{S}}^c_T = {\widehat{S}}^c_{T-h_{\ell-1}} + a({\widehat{S}}^c_{T-h_{\ell-1}})\, h_{\ell-1} +  b({\widehat{S}}^c_{T-h_{\ell-1}})\, \left(\Delta W_{N-1} {+} \Delta W_N\right),
\]
from which we obtain
\begin{eqnarray*}
{\widehat{P}}^c_{\ell-1}
&=& \bbE\left[ H({\widehat{S}}^c_T {-} K)\ |\ {\widehat{S}}^c_{T-h_{\ell-1}},\Delta W_{N-1}\right] \\
&=& \Phi\left(\frac{{\widehat{S}}^c_{T-h_{\ell-1}} + a({\widehat{S}}^c_{T-h_{\ell-1}})\, h_{\ell-1} +  b({\widehat{S}}^c_{T-h_{\ell-1}})\, \Delta W_{N-1}  - K}{b({\widehat{S}}^c_{T-h_{\ell-1}})\ \sqrt{h_\ell}}\right).
\end{eqnarray*}

With ${\widehat{Y}}_\ell \equiv {\widehat{P}}_\ell{-} {\widehat{P}}_{\ell-1}$, numerical analysis \cite{Gi_gdr19}
proves that $V_\ell \approx O(h^{3/2}_\ell)$ so the MLMC complexity is
$O(\varepsilon^{-2})$.
Heuristically, this is because there is an $O(h_\ell^{1/2})$ probability
of paths being within $O(h_\ell^{1/2})$ of the strike $K$, and for these
\[
{\widehat{S}}^f_{T-h_{\ell-1}} - {\widehat{S}}^c_{T-h_{\ell-1}} = O(h_\ell), ~~~
\frac{\partial {\widehat{P}}}{\partial {\widehat{S}}} = O(h_\ell^{-1/2})
~~~\Longrightarrow ~~~
{\widehat{P}}_\ell - {\widehat{P}}_{\ell-1} = O(h_\ell^{1/2}),
\]
so $V_\ell \approx O( h_\ell^{1/2} \times (h_\ell^{1/2})^2) = O(h_\ell^{3/2})$.
In addition, the kurtosis is improved to $O(h_\ell^{-1/2})$.
Unfortunately, this approach does not help when the Euler-Maruyama discretisation
is used for the entire path since
${\widehat{S}}^f_{T-h_{\ell-1}} - {\widehat{S}}^c_{T-h_{\ell-1}} = O(h_\ell^{1/2})$
and so ${\widehat{P}}_\ell - {\widehat{P}}_{\ell-1} = O(1)$

The use of this kind of conditional expectation is a standard technique
for smoothing the payoff to enable IPA/pathwise sensitivity calculations
\cite{Gi_glasserman04}.
Another example is a down-and-out barrier option, where the option is
knocked out if the path drops below a certain value. In this case the
payoff can be smoothed by computing the probability of this happening,
conditional on the computed path approximations at discrete timesteps
\cite{Gi_glasserman04}.
Again, this works well for MLMC when using the first order Milstein
discretisation \cite{Gi_giles08b,Gi_gdr19}, but it does not help with the
Euler-Maruyama discretisation.

A different kind of conditional expectation smoothing was introduced
by Achtsis, Cools \& Nuyens \cite{Gi_acn13} and
Bayer, Siebenmorgen \& Tempone \cite{Gi_bst18} to improve the convergence
of QMC computations, and then used by Bayer, Ben Hammouda \& Tempone
\cite{Gi_bht20} to improve the MLMC variance for digital options.

In its simplest form, they split the random inputs for the numerical
simulation into a scalar $Z$ and the remainder $Z_r$, and express
the desired MLMC level $\ell$ expectation as
\[
\bbE[ {\widehat{P}}_\ell {-} {\widehat{P}}_{\ell-1} ] = 
\bbE\left[\ \bbE[ {\widehat{P}}_\ell {-} {\widehat{P}}_{\ell-1} \,|\, Z_r ] \ \right]
\]
and observe that in many financial applications it is possible to
perform this split in a way such that the conditional expectations
$
\bbE[ {\widehat{P}}_\ell \,|\, Z_r ], 
\bbE[ {\widehat{P}}_{\ell-1} \,|\, Z_r ]
$
are smooth functions of $Z_r$, and can be evaluated analytically
or very accurately by 1D numerical quadrature when there is just
a single discontinuity with respect to changes in $Z$.

For a scalar SDE, $Z$ could be the terminal value of the driving
Brownian motion, in which case $Z_r$ would represent the other
Normal random variables required for a Brownian Bridge construction
of the Brownian increments.

\section{Change of measure}

Another approach to treating digital options using the Milstein
discretisation is to use a change of measure \cite{Gi_burgos14,Gi_giles15},
which has connections to the Likelihood Ratio Method (LRM) that
is used for sensitivity analysis \cite{Gi_ecuyer95}.

\begin{figure}[b]
\begin{center}
\includegraphics[width=0.8\textwidth]{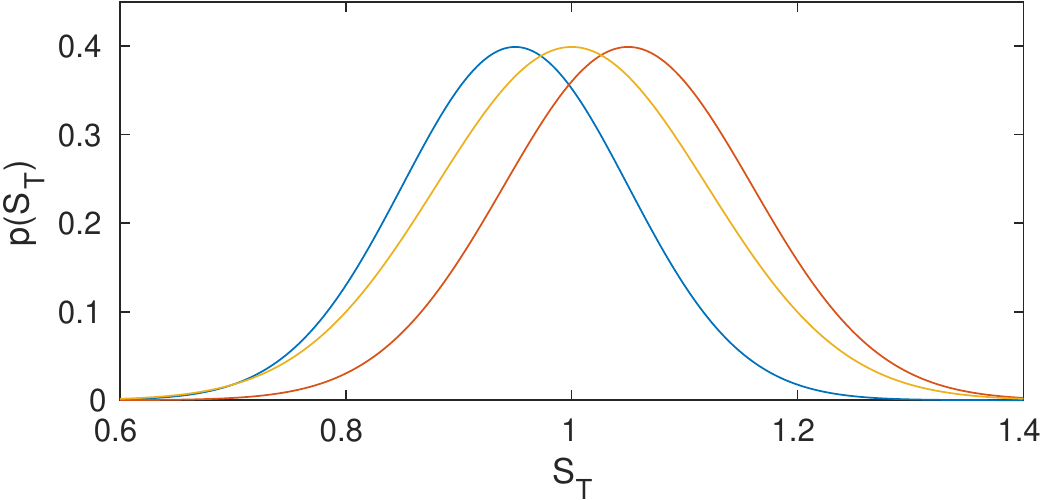}
\end{center}
\caption{Coarse and fine path conditional Gaussian distributions,
plus third common distribution}
\label{Gi_fig:measures}
\end{figure}

For both the fine and coarse paths, we have conditional Gaussian
distributions for ${\widehat{S}}_T$, with slightly different means and variances,
as illustrated in Figure \ref{Gi_fig:measures}.
We can therefore perform a change of measure to the same Gaussian 
distribution with mean $\mu$ and variance $\sigma^2$, also illustrated
in Figure \ref{Gi_fig:measures}, and then pick the same sample ${\widehat{S}}_T$
for both paths from this common Gaussian distribution.

The resulting MLMC estimator is
\[
  {\widehat{Y}}_\ell = f({\widehat{S}}_T)\ (R_\ell - R_{\ell-1})
\]
where $R_\ell, R_{\ell-1}$ are the respective Radon-Nikodym derivatives for
the fine and coarse paths. For the scalar SDE (\ref{Gi_eq:SDE}), $R_\ell$ is
\[
R_\ell = \frac{\sigma}{b({\widehat{S}}^f_{T-h_\ell}) \sqrt{h_\ell}}\
\frac{\exp\left(-\,({\widehat{S}}_T - {\widehat{S}}^f_{T-h_\ell} - a({\widehat{S}}^f_{T-h_\ell})\, h_\ell )^2
  \, /\, (2\,b^2({\widehat{S}}^f_{T-h_\ell})\, h_\ell) \right)}
     {\exp\left(-\,({\widehat{S}}_T{-}\mu)^2 \,/\, (2\sigma^2)\right)}
\]
and $R_{\ell-1}$ is defined similarly.  It can be shown that the difference
$R_\ell - R_{\ell-1}$ is approximately $O(h_\ell^{1/2})$, which implies that
$V_\ell \approx O(h_\ell)$.  To improve the variance we note that the conditional
expected value of Radon-Nikodym derivatives is always 1,
i.e.~$\bbE[R_\ell \ | \ {\widehat{S}}^f_{T-h_\ell}]
= \bbE[R_{\ell-1} \ | \ {\widehat{S}}^c_{T-h_{\ell-1}}, \Delta W_{N-1}] = 1$,
and therefore we can change the definition of ${\widehat{Y}}_\ell$ to
\[
{\widehat{Y}}_\ell \ =\ \left( f({\widehat{S}}_T) - f(\mu) \right) (R_\ell - R_{\ell-1})
\]
without changing its expected value.  This estimator is now non-zero
only when ${\widehat{S}}_T$ and $\mu$ are on opposite sides of the strike $K$,
which occurs for an $O(h_\ell^{1/2})$ fraction of coarse/fine paths.
Hence the new MLMC variance $V_\ell$ is approximately $O(h^{3/2}_\ell)$,
as with the use of the analytic conditional expectation.

The benefit of this approach is that it works well in multiple
dimensions when it is often not possible to evaluate the analytic
conditional expectation \cite{Gi_burgos14,Gi_giles15}.  However, again
it does not help with the full path Euler-Maruyama discretisation
because that gives $R_\ell - R_{\ell-1} = O(1)$.

\begin{figure}[t]
\begin{center}
\includegraphics[width=0.8\textwidth]{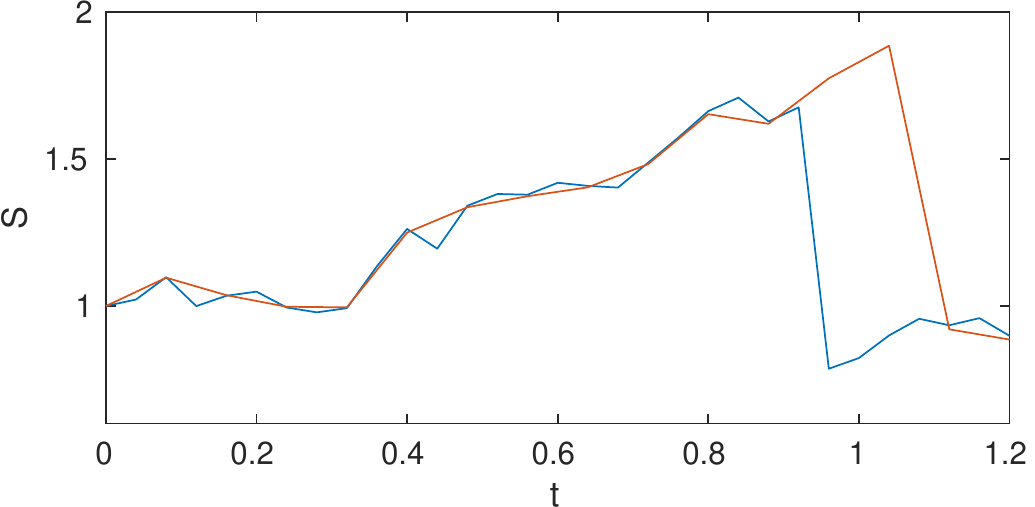}
\end{center}
\caption{Illustration of coarse and fine jump-diffusion paths
         with jumps before and after $T=1$.}
\label{Gi_fig:jumps}
\end{figure}

An earlier use of a change of measure in an MLMC computation
was by Xia \cite{Gi_xia14,Gi_xg12} for a Merton-style jump-diffusion 
SDE with a path-dependent jump rate $\lambda(S,t)$. The challenge
in this application, as illustrated in Figure \ref{Gi_fig:jumps},
is that the coarse and fine paths will jump at different times,
and one might jump just before the final time $T$, and the other
just after, leading to a large difference in the computed value of
$f(S_T)$.  The path-dependent jump rate was treated by using
the thinning technique of Glasserman \& Merener \cite{Gi_gm04},
over-sampling possible jump times using a uniform rate
$\lambda_{sup} > \lambda(S,t)$ and then using an acceptance/rejection
step to select the real jump times.  Xia modified this with
a change of measure to ensure the same acceptance/rejection
decision for both the fine and coarse paths so that they both
jump at the same time. This leads to an estimator of the form
\[
{\widehat{Y}}_\ell = {\widehat{P}}_\ell\, R_\ell - {\widehat{P}}_{\ell-1}\, R_{\ell-1}.
\]
When combined with a first order Milstein discretisation of
the SDE between the jump times, this gives $V_\ell = O(h_\ell^2)$
for Lipschitz payoff functions such as a standard put or call
option \cite{Gi_xia14,Gi_xg12}.

\section{Splitting}

Returning to the challenge of digital options arising from the
solution of an SDE, a third approach is to use path-splitting
to generate an unbiased estimate of the conditional expectation
introduced in Section 4 \cite{Gi_giles15}.

\begin{figure}[b]
\begin{center}
\includegraphics[width=0.8\textwidth]{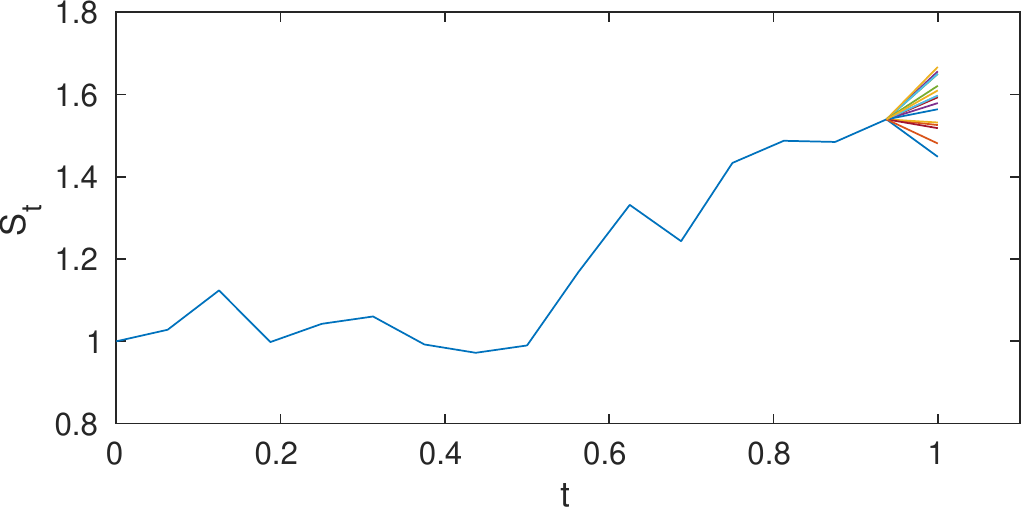}
\end{center}
\caption{Path splitting in final timestep to estimate conditional expectation}
\label{Gi_fig:splitting}
\end{figure}

This is a variant of the general splitting technique \cite{Gi_ag07}.
As illustrated in Figure \ref{Gi_fig:splitting}, it involves
performing a standard fine path simulation up until one timestep
before the final time $T$, and then performing multiple independent
simulations of the final timestep, averaging the payoff for each
of these to get an approximation of the conditional expectation.
The same is done for the coarse path except that each of the
splits uses the same $\Delta W_{N-1}$ that was used for the second
to last fine path timestep.

Since the computational cost of the path up to the splitting time
is $O(h_\ell^{-1})$, it means that up to $O(h_\ell^{-1})$ splits can
be used without increasing the path cost significantly.
If $M_\ell$ splits are used, then the standard splitting variance
analysis gives
\[
\bbV[{\widehat{Y}}_\ell] =
\bbV\left[ \bbE[ {\widehat{P}}_\ell {-}{\widehat{P}}_{\ell-1} \ | \ \{\Delta W_n\}_{n<N}] \right]
 + M_\ell^{-1} 
\bbE\left[ \bbV[ {\widehat{P}}_\ell {-}{\widehat{P}}_{\ell-1} \ | \ \{\Delta W_n\}_{n<N}] \right].
\]
As discussed previously
$\bbV\left[ \bbE[ {\widehat{P}}_\ell {-}{\widehat{P}}_{\ell-1} \ | \ \{\Delta W_n\}_{n<N}] \right]
= O(h_\ell^{3/2})$, and similarly it can be argued that
$\bbE\left[ \bbV[ {\widehat{P}}_\ell {-}{\widehat{P}}_{\ell-1} \ | \ \{\Delta W_n\}_{n<N}] \right]
= O(h_\ell)$.
Therefore choosing $M_\ell$ to lie between $O(h_\ell^{-1})$ and
$O(h_\ell^{-1/2})$ ensures the benefits of the splitting are obtained
without significantly increasing the computational cost per sample.

As an additional bonus, one can use the more accurate Milstein
discretisation for the final timestep, instead of switching to
the Euler-Maruyama discretisation.
Burgos \cite{Gi_burgos14,Gi_bg12} gives more details of the analysis,
and also used the same approach for pathwise sensitivity analysis
for a variety of financial options.

Giles \& Bernal \cite{Gi_gb18} also used splitting for Feynman-Kac
functionals arising for stopped diffusions, SDE simulations which
terminate when the solution path leaves a prescribed domain.
The issue here is that when a fine path exits, there is 
an $O(h_\ell^{1/2})$ probability that the corresponding coarse 
path does not leave until much later. This is addressed 
by estimating a conditional expectation by splitting the
coarse path into $O(h_\ell^{-1/2})$ independent sub-simulations
which continue until each of them leaves the domain.
$V_\ell$ is improved from $O(h^{1/2}_\ell)$ to approximately
$O(h_\ell)$ without a significant increase in the cost per
sample, and finally the MLMC complexity achieved is
$O(\varepsilon^{-2}|\log\varepsilon|^3)$.

None of the three methods introduced so far (conditional
expectation, change of measure, splitting) helps when
using the Euler-Maruyama discretisation.  For this, a new
method has recently been developed by Giles \& Haji-Ali
\cite{Gi_gh22b}.

\begin{figure}[b]
\begin{center}
\includegraphics[width=0.8\textwidth]{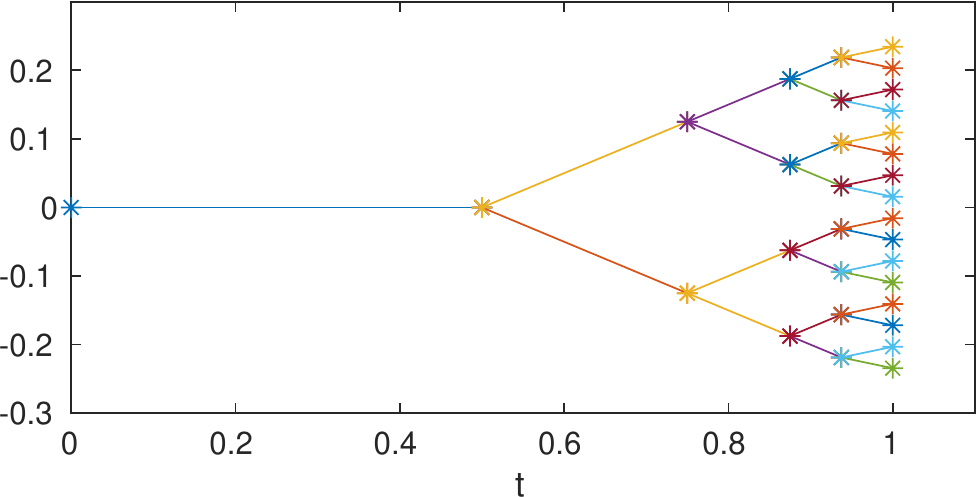}
\end{center}
\caption{Repeated path splitting to estimate conditional expectation}
\label{Gi_fig:new-splitting}
\end{figure}

It again uses splitting, but inspired by
the simulation of branching diffusions, it considers
splits at multiple deterministic times, as illustrated
in Figure \ref{Gi_fig:new-splitting} which shows the logical
structure of a set of split paths.  Here we are considering
a simulation on the unit time interval.  A single pair of
fine/coarse paths is calculated up to time $t=1/2$, with
the number of fine timesteps being $\frac{1}{2}h_\ell^{-1}$.
This simulation is then split into two separate independent
simulations up to time $t=3/4$, with the two simulations
between them accounting for an additional $\frac{1}{2}h_\ell^{-1}$
fine timesteps. There are further splits at $t=3/4$, then at
$t=7/8$, and so on, with the final split when there is just
one coarse timestep left.

The total number of fine timesteps simulated is
$O(h_\ell^{-1}|\log h_\ell|)$ so the computational cost is only
slightly increased compared to the original method with a
single pair of fine/coarse paths.  ${\widehat{Y}}_\ell$ is defined to be
the average of the values ${\widehat{P}}_\ell{-}{\widehat{P}}_{\ell-1}$ for each of
the final paths, and it can be proved that its variance is
$O(h_\ell)$, the same asymptotic order of convergence as
for Lipschitz payoff functions \cite{Gi_gh22b}.  The kurtosis
is also improved, so this technique fully addresses the
challenge of using MLMC with the Euler-Maruyama discretisation
to estimate digital option values.

\section{Adaptive sampling}

We return now to the challenge of estimating the nested
expectation
$\displaystyle
\bbE \left[\,  H\left(\bbE[Z | X]\right)  \, \right] 
$
and we note that we only need an accurate estimate of the
inner conditional expectation $\bbE[Z | X]$ when it is near
zero.
This observation is the basis for the development of adaptive 
sampling by Broadie, Du \& Moallemi \cite{Gi_bdm11} within a
standard Monte Carlo procedure.  This was then extended to
adaptive sampling combined with MLMC by Giles \& Haji-Ali \cite{Gi_gh19}
by defining the number of inner samples $M_\ell$ on level $\ell$ to be
\begin{itemize}
\item  $M_\ell=2^\ell M_0$ inner samples when $|\bbE[Z | X]| \gg \sqrt{\bbV[Z|X]/(2^\ell M_0)}$

  \vspace{0.05in}
  
  This is the smallest number of samples used on level $\ell$.
  $\sqrt{\bbV[Z|X]/(2^\ell M_0)}$ is the standard deviation of the 
  Monte Carlo estimate for $\bbE[Z | X]$, so the inequality means 
  that this number of samples is sufficient to be very sure that 
  the estimate has the correct sign.

  \vspace{0.05in}

\item  $M_\ell=4^\ell M_0$ inner samples when $|\bbE[Z | X]| = O(\sqrt{\bbV[Z|X]/(4^\ell M_0)})$

  \vspace{0.05in}
  
  This is the maximum number of samples used on level $\ell$.  In this case,
  the estimate of $\bbE[Z | X]$ may have the incorrect sign, but this will
  only happen when $|\bbE[Z | X]| = O(2^{-\ell})$ which occurs with probability
  $O(2^{-\ell})$.  Likewise, the total cost of the higher number of
  samples in this region is $O(2^{-\ell} \times 4^\ell) = O(2^\ell)$, so it
  does not significantly increase the overall average cost.

  \vspace{0.05in}
  
\item  $2^\ell M_0<M_\ell<4^\ell M_0$ for intermediate values

  \vspace{0.05in}

  In this region the number of samples is chosen to be very sure
  that the estimate of $\bbE[Z | X]$ has the correct sign, and at the same
  time the total cost is $O(2^\ell)$.
  
\end{itemize}

\begin{figure}[b]
\begin{center}
\includegraphics[width=0.8\textwidth]{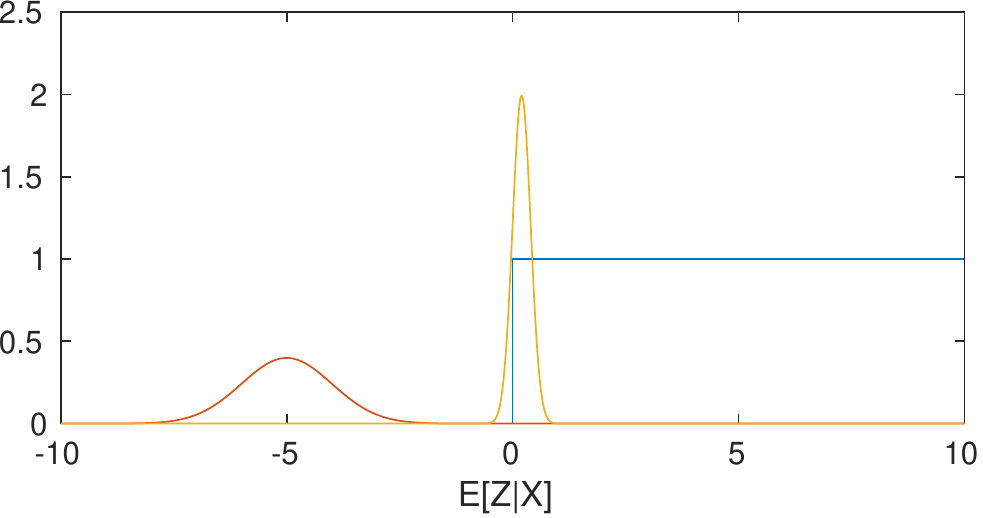}
\end{center}
\caption{Error distributions for two conditional expectations with
  i) few samples being needed to ensure the correct sign (left), and
  ii) many samples being insufficient to ensure the correct sign (centre).
The blue line represents the Heaviside step function.}
\end{figure}

Overall, this adaptive sampling approach leads to $C_\ell\sim 2^\ell$,
$V_\ell\sim 2^{-\ell}$ and hence a complexity of roughly $O(\varepsilon^{-2})$
\cite{Gi_gh19}.  However, the kurtosis is $O(2^\ell)$ since only an
$O(2^{-\ell})$ fraction of the outer samples give non-zero values
for ${\widehat{Y}}_\ell$.

Haji-Ali, Spence \& Teckentrup \cite{Gi_hst22} have further extended
this to estimate quantities of the form
\[
\bbP[G\in \Omega] \equiv \bbE[{\bsone}_{G\in \Omega}]
\]
where $G$ is a $d$-dimensional random variable which cannot
be sampled directly. In their paper they consider in particular
the two challenges in this article. In the context of the digital
option with the Euler-Maruyama discretisation on the unit time interval,
the adaptive sampling varies the timestep used on level $\ell$ so that
\begin{itemize}
\item  $h_\ell=2^{-\ell}$ when $|{\widehat{S}}_\ell-K|$ is large compared to
       the strong error in the path approximation

  \vspace{0.05in}
  
\item  $h_\ell=4^{-\ell}$ when $|{\widehat{S}}_\ell-K|$ is of the same order
       as the strong error

  \vspace{0.05in}

\item  $2^{-\ell} < h_\ell <4^{-\ell}$ for intermediate values
\end{itemize}
A Brownian bridge construction is used when the timestep
needs to be refined as part of the adaptation procedure from
its initial value $h_\ell=2^{-\ell}$.
The adaptation again leads to $C_\ell\sim 2^\ell$,
$V_\ell\sim 2^{-\ell}$ and hence a complexity of roughly
$O(\varepsilon^{-2})$, but there is again a high kurtosis \cite{Gi_hst22}.

In earlier research, Elfverson, Hellman \& Malqvist \cite{Gi_ehm16}
considered estimation of $\bbE[H(X)]$ where $X$ cannot be sampled
exactly but there is a sequence of approximations $X'_0, X'_1,
X'_2, \ldots X$ of increasing accuracy and increasing cost.
Motivated by PDE applications with a well-behaved truncation error
so that there are uniform geometric bounds on $|X'_j-X|$, level
$\ell$ in their method uses
\[
{\widehat{X}}_\ell = X'_j, ~~~ j = \min\{ \ell, \min j: |{\widehat{X}}'_j-X| < |X| \}
\]
and achieves similarly good MLMC benefits. This idea is essentially
the same as in the work of Haji-Ali {\it et al} but requiring a
uniform bound on $|X'_j{-}X|$ is significantly more restrictive than
the bounds on $\bbE[\, |X'_j{-}X|^q]$ for some $q{>}2$ required by
Haji-Ali {\it et al}.

A final comment is that the  analysis of Haji-Ali, Spence \&
Teckentrup can be generalised to a product of an indicator function
and a Lipschitz function,  $\bbE[{\bsone}_{G\in \Omega}f(S)]$,
and so can handle barrier options.  Furthermore, 
Haji-Ali \& Spence have extended the adaptive sampling methodology
to an extremely challenging triply-nested expectation which arises
in mathematical finance \cite{Gi_ghs23}. By incorporating the randomised
MLMC treatment of Rhee \& Glynn \cite{Gi_rg15} to handle the time
discretisation of the underlying SDEs as well as the sampling for
the inner conditional expectations, they achieve an overall complexity
of approximately $O(\varepsilon^{-2})$ which is very impressive for such a
difficult application.

\section{Conclusions}

It is worth repeating that in most MLMC applications the output
quantity of interest is a Lipschitz function of the intermediate
simulation quantities, so good strong convergence for the intermediate
quantities leads automatically to a good rate of convergence of
the MLMC variance $V_\ell$.

For those applications in which the function is discontinuous,
this article shows there is an extensive literature with a
variety of different approaches to improve the MLMC variance
and try to recover the optimal $O(\varepsilon^{-2})$ complexity. It
is notable that many of these methods have adapted ideas from
Monte Carlo sensitivity analysis which also has problems with
discontinuous functionals.  It is hoped that this survey will
assist future researchers facing similar challenges in other
new application areas.

\begin{acknowledgement}
  This paper is based on research with many students, postdocs
  and other collaborators and I am grateful to all of them.
  Funding for the research has been provided by the UK
  Engineering and Physical Sciences Research Council
  through grants EP/E031455/1, EP/H05183X/1 and EP/P020720/2
  as well as the Hong Kong Innovation and Technology Commission
  (InnoHK Project CIMDA).
  The paper was written while visiting the Oden Institute for
  Computational Engineering \& Sciences at UT Austin, and I
  thank my hosts for their warm hospitality.
\end{acknowledgement}
%



\end{document}